\documentclass[]{amsart}

\usepackage{amssymb}
\usepackage{amsmath}
\usepackage{indentfirst}
\usepackage{amscd}
\usepackage{graphicx}

\newtheorem{theorem}{Theorem}[section]
\newtheorem{lem}[theorem]{Lemma}

\newtheorem{cor}[theorem]{Corollary}

\newcommand{\R}{\ensuremath{\mathbb{R}}}

\newcommand{\Z}{\ensuremath{\mathbb{Z}}}

\newcommand{\p}{\ensuremath{\mathbf{v_+}}}
\newcommand{\m}{\ensuremath{\mathbf{v_-}}}

\begin{document}


\subjclass[2000]{Primary 57Q45; Secondary 57M25.}
\date{First edition: Feb. 8, 2005.  This edition: Feb. 27, 2005.}
\keywords{Khovanov cohomology, surface-knot, Khovanov-Jacobsson number.}


\title[Khovanov-Jacobsson numbers and invariants of surface-knots]
{Khovanov-Jacobsson numbers and invariants of surface-knots derived from
Bar-Natan's theory}


\author{Kokoro Tanaka}
\address{Graduate School of Mathematical Sciences, 
University of Tokyo, 3-8-1 Komaba Meguro, 
Tokyo 153-8914, Japan}
\email{k-tanaka@ms.u-tokyo.ac.jp}



\begin{abstract}
Khovanov introduced a cohomology 
theory for oriented classical links whose graded 
Euler characteristic is the Jones polynomial. 
Since Khovanov's theory is functorial
for link cobordisms between classical links, 
we obtain an invariant of a surface-knot, called the
{\it Khovanov-Jacobsson number}, by considering
the surface-knot as a link cobordism between empty links.
In this paper, we define an invariant of a surface-knot
which is a generalization of the Khovanov-Jacobsson number 
by using Bar-Natan's theory, and prove that 
any $T^2$-knot has the trivial Khovanov-Jacobsson number.
\end{abstract}

\maketitle


\section{Introduction}
Khovanov \cite{Kh1} introduced a cohomology 
theory for oriented classical links which values in 
graded $\Z$-modules 
and whose graded Euler characteristic is the Jones polynomial.
We denote Khovanov's cohomology groups of an oriented
link $L$ by $H(L,\mathcal{F})$ $(=\bigoplus H^i (L,\mathcal{F}) )$. 
His theory is powerful for classical links; for example,
Bar-Natan \cite{BN1} and  Wehrli \cite{W}
showed that Khovanov's cohomology is stronger than 
the Jones polynomial, and Rasmussen \cite{Ras} gave a 
combinatorial proof of the Milnor conjecture
by using a variant of Khovanov's theory defined by Lee \cite{Lee}.

Jacobsson \cite{Jac} and Khovanov \cite{Kh2} proved that
Khovanov's theory is functorial 
for link cobordisms in the following sense: 
a link cobordism $S \subset \R^3 \times [0,1]$
between classical links $L_0$ and $L_1$
induces a homomorphism $\phi_S : 
H(L_0,\mathcal{F}) \rightarrow H(L_1,\mathcal{F})$,
well-defined up to overall minus sign, 
under ambient isotopy of $S$ rel $\partial S$.
A {\it surface-knot} $F$ is a closed connected oriented
surface embedded locally flatly in $\R^4$,
and can be considered as a link cobordism 
between empty links. 
Then the induced map $\phi_F :
H(\emptyset,\mathcal{F} ) \rightarrow 
H(\emptyset,\mathcal{F} )$, up to overall minus sign,
gives an invariant of the surface-knot $F$.
Since the cohomology group $H(\emptyset, \mathcal{F} )$
of empty link $\emptyset$ is $\Z$, 
the map $\phi_F$ is an endomorphism of $\Z$.
Hence we obtain an invariant of the surface-knot $F$ 
defined as $\mid \phi_F (1) \mid \in \Z,$
and denote it by $KJ(F)$. 
This invariant is called the 
{\it Khovanov-Jacobsson number} in \cite{CSS}.

As far as the author knows,
there are a few result on the computation of 
the Khovanov-Jacobsson numbers of surface-knots; 
see \cite{CSS}, for example.
It follows from a simple observation that
$KJ(F)=0$ for any surface-knot $F$ with $\chi(F) \neq 0$
and that $KJ(F)=2$ for a trivial $T^2$-knot $F$
(a surface-knot with $\chi(F)=0$).
It seems to be hard to compute the Khovanov-Jacobsson number
in general, but Carter, Saito and Satoh \cite{CSS} proved that 
$KJ(F)=2$ for any $T^2$-knot $F$ obtained 
from a spun/twist-spun $S^2$-knot by attaching a $1$-handle.
They also proved that $KJ(F)=2$ for any 
pseudo-ribbon $T^2$-knot $F$.

In this paper, we define an invariant 
$BN(F) \in \Z[t]$ of a surface-knot $F$
by using a variant of Khovanov's theory
defined by Bar-Natan \cite{BN2}. 
This invariant is a generalization of 
the Khovanov-Jacobsson number such that
$$BN(F) |_{t=0} = KJ(F).$$
The main result of this paper is that 
the invariant $BN(F)$ is trivial for
any surface-knot $F$, and hence it turns out that
the Khovanov-Jacobsson number is trivial 
for any $T^2$-knot\footnote{The author has subsequently learned that 
Jacob Rasmussen \cite{Ras-KJ} has a different proof of Corollary~\ref{cor-KJ}
using Lee's theory.} :

\begin{theorem}\label{th-main}
For any surface-knot $F$ of genus $g$ $(g\geq 0)$, 
we have the following.
\begin{itemize}
\item[(i)]
If $g$ is an even integer, then we have $BN(F)=0$.
\item[(ii)]
If $g$ is an odd integer, then we have $BN(F)=2^g t^{(g-1)/2}$. 
\end{itemize}
\end{theorem}

\begin{cor}\label{cor-KJ}
For any $T^2$-knot, we have $KJ(F)=2$. 
\end{cor}

This paper is organized as follows.
In Section~\ref{sec-Kh} and ~\ref{sec-KJ},
we review Khovanov's cohomology theory for oriented classical links
and the Khovanov-Jaconsson numbers of surface-knots respectively.
In Section~\ref{sec-BN}, we define the surface-knot
invariant $BN(F)$.
Section~\ref{sec-pf} is devoted to the proof of 
Theorem~\ref{th-main}.

\section{Khovanov's cohomology theory}\label{sec-Kh}
In this section, we briefly recall Khovanov's 
cohomology theory for oriented classical links \cite{Kh1}.
See also \cite{BN1}, \cite{BN2}, \cite[Section $2$]{Lee},
\cite[Section $2$]{Ras}, for example.

\subsection{Graded $\Z$-module $V$ and TQFT $\mathcal{F}$}
Let $V$ be a free graded $\Z$-module of rank two generated by
$\p$ and $\m$ with 
$$\text{deg}(\p)=1 \ \, \text{and} \ \, \text{deg}(\m)=-1 .$$
We give the graded $\Z$-module $V$ 
a Frobenius algebra structure with
a multiplication $m$, a comultiplication $\Delta$,
a unit $\iota$, and a counit $\epsilon$ defined by
$$\begin{array}{lll}
m(\p \otimes \p) = \p & \quad &
\Delta(\p) = \p \otimes \m + \m \otimes \p  \\
m(\p \otimes \m) = m(\m \otimes \p) = \m & \quad &
\Delta(\m) = \m \otimes \m \\
m(\m \otimes \m) = 0 & & \\
\iota(1) = \p &  & \epsilon(\p) = 0 \ \ \   \epsilon(\m) = 1 .
\end{array}$$
The structure maps $m$, $\Delta$, $\iota$ and $\epsilon$ are
graded maps of degree $-1$, $-1$, $1$ and $1$ respectively. 

Khovanov's cohomology theory is based on a
$(1+1)$-dimensional TQFT $\mathcal{F}$,
a monoidal functor from oriented $(1+1)$-cobordisms to 
graded $\Z$-modules, associated to $V$.
The Frobenius algebra $V$ defines $\mathcal{F}$ by assigning
$\Z$ to an empty $1$-manifold,
$V$ to a single circle, $V \otimes V$ to a disjoint union 
of two circles, and so on.
The structure maps are assigned to elementary cobordisms 
such that
$$\begin{array}{ll}
\mathcal{F}\left(
 \begin{minipage}{15pt}
  \includegraphics[width=15pt]{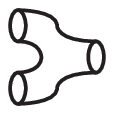}
 \end{minipage}\right)=
m : V \otimes V \rightarrow V, &
\mathcal{F}\left(
 \begin{minipage}{15pt}
  \includegraphics[width=15pt]{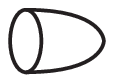}
 \end{minipage}\right)=
\iota : \Z \rightarrow V, \\
\mathcal{F}\left(
 \begin{minipage}{15pt}
  \includegraphics[width=15pt]{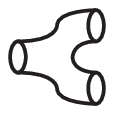}
 \end{minipage}\right)=
\Delta : V \rightarrow V \otimes V, &
\mathcal{F}\left(
 \begin{minipage}{15pt}
  \includegraphics[width=15pt]{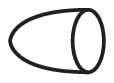}
 \end{minipage}\right)=
\epsilon : V \rightarrow \Z \ .\end{array}$$

\subsection{Cube of resolutions}
Let $L$ be an oriented link, and 
$D$ an oriented link diagram of $L$ with $n$ 
crossings labeled by $1, 2, \dots, n$.
A double point of $D$ can be resolved in two ways:
one is the $0$-smoothing 
\begin{minipage}[10pt]{10pt}\begin{picture}(10,10)
 \qbezier(1,1)(1,1)(9,9)
 \qbezier(1,9)(1,9)(3,7)
 \qbezier(7,3)(7,3)(9,1)
\end{picture}\end{minipage}
$\rightsquigarrow$
\begin{minipage}[10pt]{10pt}\begin{picture}(10,10)
 \qbezier(1,1)(7,5)(1,9)
 \qbezier(9,1)(3,5)(9,9)
\end{picture}\end{minipage},
and the other is the $1$-smoothing
\begin{minipage}[10pt]{10pt}\begin{picture}(10,10)
 \qbezier(1,1)(1,1)(9,9)
 \qbezier(1,9)(1,9)(3,7)
 \qbezier(7,3)(7,3)(9,1)
\end{picture}\end{minipage}
$\rightsquigarrow$
\begin{minipage}[10pt]{10pt}\begin{picture}(10,10)
 \qbezier(1,1)(5,7)(9,1)
 \qbezier(1,9)(5,3)(9,9)
\end{picture}\end{minipage}.
We construct an $n$-dimensional cube $[0,1]^n$,
called the {\it cube of resolutions}, 
from $D$ such that each vertex $v$
is decorated by a complete resolution of $D$ and 
each edge $e$ is decorated by a cobordism.

To each vertex $v$ of the cube, 
we associate the complete resolution $D_v$ as follows:
the $i$th crossing is resolved by the $0$-smoothing 
if the $i$th coordinate of $v$ is $0$,
and by the $1$-smoothing if it is $1$.
Then $D_v$ is a collection of simple closed curves.

Each edge $e$ of the cube  
can be represented by sequences in $\{ 0,1, \star \}^n$
with just one $\star$, and has the two end vertices
$v_e(0)$ and $v_e(1)$: the vertex $v_e(0)$
is obtained by substituting $0$ to $\star$,
and the vertex $v_e(1)$ is obtained 
by substituting $1$ to $\star$.
To an edge $e$ for which the $j$th 
coordinate is $\star$, we associate the cobordism 
$S_e : D_{v_e(0)} \rightarrow D_{v_e(1)}$ as follows:
we remove a neighborhood of the $j$th crossing, 
assign a product cobordism, and fill the saddle cobordism
between the $0$- and $1$-smoothings around the $j$th crossing.
The cobordism $S_e$ is either of the following two types:
(i) two circles of $D_{v_e(0)}$ merge into
one circle of $D_{v_e(1)}$, or
(ii) one circle of $D_{v_e(0)}$
splits into two circles of $D_{v_e(1)}$.

\subsection{Cube of modules}
Applying $\mathcal{F}$ to the cube of resolutions,
we construct another $n$-dimensional cube $[0,1]^n$,
called the {\it cube of modules}.
Each vertex $v$ is replaced by a graded $\Z$-module
$\mathcal{F}(D_v)$ and each edge $e$ is replaced by 
a homomorphism 
$$\mathcal{F}(S_e) : 
\mathcal{F}\left( D_{v_e(0)}\right)
\rightarrow \mathcal{F}\left( D_{v_e(1)}\right).$$ 
The homomorphism $\mathcal{F}(S_e)$ is induced by a map
$m:V \otimes V  \rightarrow V$
if the cobordism $S_e$ is of type (i), and 
is induced by a map $\Delta : V \rightarrow V \otimes V$
if it is of type (ii).

\subsection{Khovanov's cohomology groups}
For a graded $\Z$-module $M$ and an integer $k$,
let $M \{k\}$ denote the graded $\Z$-module
of which the $j$th graded component is the $(j-k)$th
graded component of $M$.
We note that $M\{k\}$ is identical to $M$ as $\Z$-modules.
Khovanov's cochain complex $C(D,\mathcal{F})$
of graded $\Z$-modules is obtained from the cube of modules 
as follows. The underlying group of $C(D,\mathcal{F})$ 
$(=\oplus \, C^i(D,\mathcal{F}) )$ is defined by
$$\begin{array}{ll}
C^i(D,\mathcal{F}) &=\bigoplus\limits_{v : \ i=|v| - n_-}
\mathcal{F}(D_v) \bigl\{ (|v|-n_-) + (n_+ - n_-) \bigr\} ,
\end{array}$$
where $|v|$ is the sum of all coordinates of $v$
and $n_+$ (resp. $n_-$) is the number of 
positive (resp. negative) crossings of $D$.
The coboundary map $d$ for an element $x$ of $\mathcal{F}(D_v)$
is defined by
$$d(x) =\sum\limits_{e_i : \ v_{e_i}(0) = v}
(-1)^{w(e_i)} \mathcal{F}(S_{e_i}) (x),$$
where $w(e)$ is the sum of all coordinates of an
edge $e$ after $\star$. 
We give some remarks about the coboundary
map $d$ in Section~\ref{subsec-remark}.

The cohomology groups $H^i(D,\mathcal{F})$ of the 
cochain complex $C(D,\mathcal{F})$ are graded $\Z$-modules.
Khovanov \cite{Kh1} proved that $C(D,\mathcal{F})$ and 
$C(D',\mathcal{F})$ are cochain homotopic to each other
for any other diagram $D'$ of $L$.
To prove it, he constructed cochain maps 
$$\begin{array}{ll}
f_1 : C \left(
 \begin{minipage}{15pt}
  \includegraphics[width=15pt]{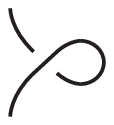}
 \end{minipage} , 
\mathcal{F} \right)
\rightarrow  C \left( 
\begin{minipage}{15pt}
  \includegraphics[width=15pt]{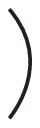}
 \end{minipage}, 
\mathcal{F} \right), &
g_1 : C \left( 
\begin{minipage}{15pt}
  \includegraphics[width=15pt]{r1-line.eps}
 \end{minipage}, 
\mathcal{F} \right)
\rightarrow  C \left( 
\begin{minipage}{15pt}
  \includegraphics[width=15pt]{r1.eps}
 \end{minipage}, 
\mathcal{F} \right), \\
f_2 : C \left( 
\begin{minipage}{15pt}
  \includegraphics[width=15pt]{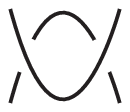}
 \end{minipage}, 
\mathcal{F} \right)
\rightarrow  C \left( 
 \begin{minipage}{15pt}
  \includegraphics[width=15pt]{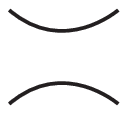}
 \end{minipage}, 
\mathcal{F} \right), &
g_2 : C \left( 
 \begin{minipage}{15pt}
  \includegraphics[width=15pt]{r2-01.eps}
 \end{minipage}, 
\mathcal{F} \right)
\rightarrow  C \left( 
 \begin{minipage}{15pt}
  \includegraphics[width=15pt]{r2.eps}
 \end{minipage}, 
\mathcal{F} \right), \\
f_3 :  C \left( 
 \begin{minipage}{15pt}
  \includegraphics[width=15pt]{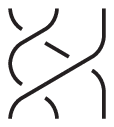}
 \end{minipage}, 
\mathcal{F} \right)
\rightarrow  C \left(\begin{minipage}{15pt}
  \includegraphics[width=15pt]{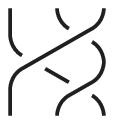}
 \end{minipage} , 
\mathcal{F} \right), &
g_3 :  C \left( 
 \begin{minipage}{15pt}
  \includegraphics[width=15pt]{r3-d.eps}
 \end{minipage}, 
\mathcal{F} \right)
\rightarrow  C \left( 
\begin{minipage}{15pt}
  \includegraphics[width=15pt]{r3-u.eps}
 \end{minipage}, 
\mathcal{F} \right),
\end{array}$$
for Reidemeister moves 
such that the maps $g_i \circ f_i$ and $f_i \circ g_i$ are 
cochain homotopic to the identities
for each $i$ $(i=1,2,3)$,
but we omit the precise definition of these maps.
(Refer to \cite{Kh1,BN1,BN2} for more details.)
This implies that the isomorphism classes of the cohomology groups 
$H^i(D,\mathcal{F})$ are invariants of $L$.
Then we denote cohomology groups of an oriented
link $L$ by $H(L,\mathcal{F})$ $(=\bigoplus H^i (L,\mathcal{F}) )$.

\subsection{Remarks}\label{subsec-remark}
We give some remarks about Khovanov's cochain complex.
\begin{itemize}
\item[(I)]
The Frobenius algebra $V$ has the following properties:
the map $m$ is associative, the map $\Delta$ is coassociative, and
$\Delta \circ m = (m \otimes \text{id}) \circ 
(\text{id} \otimes \Delta)$.
These properties can be interpreted as the commutativity
of saddle point moves in oriented $(1+1)$-cobordisms.
\item[(II)]
It follows from (I) that each two dimensional face of 
the cube of modules is commutative. 
If we replace the map $\mathcal{F}(S_{e})$ by
the map $(-1)^{w(e)} \mathcal{F}(S_{e})$ for each edge $e$,
then each two dimensional face becomes 
anti-commutative. 
\item[(III)]
The above (anti-)commutativity of the cube of modules
ensures that the map $d$ satisfies $d \circ d = 0$.
\end{itemize}

\section{Khovanov-Jacobsson numbers of surface-knots}\label{sec-KJ}

A {\it link cobordism} $S$ between classical links 
$L_0$ and $L_1$ is a compact oriented surface 
embedded properly and locally flatly in $\R^3 \times [0,1]$ 
such that $S \cap (\R^3 \times \{i\})=L_i$ for 
each $i$ $(i=0,1)$. 
In this section, we briefly recall
the map $\phi_S : H(L_0,\mathcal{F}) \rightarrow 
H(L_1,\mathcal{F})$ induced by the link cobordism $S$ 
and the Khovanov-Jaconsson numbers of surface-knots
(cf. \cite{Jac,Kh2,BN2}).

\subsection{Movie presentations}
For a fixed projection 
$\pi: \R^3 \times [0,1] \rightarrow \R^2 \times [0,1]$,
each intersection $\pi (S) \cap (\R^2 \times \{ t \})$
of a link cobordism $S$ is called a {\it still}
and denoted by $D_t$ for each $t$ $(t \in [0,1])$.
By perturbing $S$ if necessary,
we may assume that each still is a classical 
link diagram for all but finitely many critical values 
and contains at most one singular point.
When $t$ passes through a critical 
value, we see one of the following changes, called 
{\it elementary string interaction}s
(or {\it ESI}s, for short):
Reidemeister moves (R1, R2 and R3) and Morse moves
(birth, death and saddle). 
The collection $\{D_t\}_{t\in [0,1]}$ of stills
is called a {\it movie} of $S$. 
Refer to \cite{CS-book} for more details.

\subsection{The maps induced by ESIs}
We construct a map $\phi_S$ for 
a link cobordism $S$ represented by a single 
elementary string interaction.
For the $i$th Reidemeister move,
the map $\phi _S$ is defined to be ${f_i}_\ast$ or 
${g_i}_\ast$ for each $i$ $(i=1,2,3)$.
For a birth (resp. death) move,
the map $\phi _S$ is induced by $\iota$ (resp. $\epsilon$).
For a saddle move, 
the map $\phi _S$ is induced by either $m$ or $\Delta$ on 
each component of the cochain complex,
depending on whether the saddle move merges two
circle or splits one circle into two circles at 
the cube of resolutions level.
We note that the map $\phi _S$ is a graded map of degree $0$ 
for a Reidemeister move, degree $1$ for a 
birth/death move, and degree $-1$ for a saddle move.

\subsection{The maps induced by link cobordisms}
For a link cobordism $S$, 
take a movie presentation of $S$.
As mentioned above, the movie presentation of $S$ is 
represented as a collection of elementary string interactions:
$S=S_1 \cup \dots \cup S_k$.
The map $\phi _S$ is defined to be the composite 
$\phi_{S_k}\circ \dots \circ \phi_{S_1}$ of the maps
induced by elementary cobordisms $S_1, \dots ,S_k$.
Then the map $\phi _S$ is a graded map of degree $\chi (S)$, 
where $\chi(S)$ is the Euler characteristic of $S$.
Jacobsson \cite{Jac} and Khovanov \cite{Kh2} proved that
a link cobordism $S$ between classical links $L_0$ and $L_1$
induces a homomorphism 
$\phi_S : 
H(L_0,\mathcal{F}) \rightarrow H(L_1,\mathcal{F}),$
well-defined up to overall minus sign, 
under ambient isotopy of $S$ rel $\partial S$.
(Jacobsson \cite{Jac} pointed out that if we do not 
assume \lq \lq rel $\partial S$\rq \rq,
then the statement does not hold in general.)

\subsection{Khovanov-Jacobsson numbers}
A {\it surface-knot} $F$ is a closed connected oriented
surface embedded locally flatly in $\R^4$.
When we consider the surface-knot $F$ as an oriented cobordism 
between empty links, the induced map $\phi_F :
H(\emptyset,\mathcal{F} ) \rightarrow 
H(\emptyset,\mathcal{F} )$, up to overall minus sign,
gives an invariant of $F$.
Since the cohomology group $H(\emptyset, \mathcal{F} )$
of empty link $\emptyset$ is $\Z$, 
the map $\phi_F$ is an endomorphism of $\Z$.
Hence we obtain an invariant of the surface-knot $F$ 
defined as $\mid \phi_F (1) \mid \in \Z,$
and denote it by $KJ(F)$. 
This invariant is called the 
{\it Khovanov-Jacobsson number} in \cite{CSS}.
It is easy to see that
$KJ(F)=0$ for any surface-knot $F$ with $\chi(F) \neq 0$,
since the map $\phi_F$ is a graded map of degree $\chi(F)$.

\section{The surface-knot invariant derived from 
Bar-Natan's thoery}\label{sec-BN}

Bar-Natan \cite{BN2} defined several variants of Khovanov's theory.
Let $V'$ be a free graded $\Z[t]$-module of rank two generated by
$\p$ and $\m$ with 
$$\text{deg}(t)=-4, \ \, \text{deg}(\p)=1 \ \,
\text{and} \ \, \text{deg}(\m)=-1 .$$
We give $V'$ a Frobenius algebra structure with
a multiplication $m'$, a comultiplication $\Delta '$,
a unit $\iota '$, and a counit $\epsilon '$ defined by
$$\begin{array}{lll}
m'(\p \otimes \p) = \p & \quad &
\Delta '(\p) = \p \otimes \m + \m \otimes \p  \\
m'(\p \otimes \m) = m'(\m \otimes \p) = \m & \quad &
\Delta '(\m) = \m \otimes \m + t \p \otimes \p \\
m'(\m \otimes \m) = t \p .& & \\
\iota '(1) = \p  & &\epsilon '(\p)=0 \ \ \ \epsilon '(\m)=1.
\end{array}$$
One of his cohomology theories,
implicitly defined in \cite[Section 9.2]{BN2}, 
is based on a $(1+1)$-dimensional TQFT
$\mathcal{F}'$ associated to a Frobenius algebra $V'$,
and we denote cohomology groups of an oriented
link $L$ by $H(L,\mathcal{F}')$ $(=\bigoplus H^i (L,\mathcal{F}') )$.
Here the cohomology groups $H^i(D,\mathcal{F})$ are graded 
$\Z[t]$-modules. We note that the cohomology theory 
associated to $\mathcal{F}'$ is essentially the same
as that associated to the Frobenius system 
$\mathcal{F}_3$ in \cite{Kh3}, but 
the notational conventions are slightly different.

He proved that the cohomology theory 
associated to $\mathcal{F}'$
is also functorial for link cobordisms up to sign 
indeterminacy (cf. \cite{BN2}, \cite[Proposition $6$]{Kh3}).
Given a surface-knot $F$, the induced map 
$\psi_F : H(\emptyset,\mathcal{F}' ) \rightarrow 
H(\emptyset,\mathcal{F}' )$ becomes an endomorphism of $\Z[t]$.
Hence we obtain an invariant of the surface-knot $F$ 
defined as $$\mid \psi_F (1) \mid \in \Z[t],$$
and denote it by $BN(F)$. 
It follows from the definition of $V$ and $V'$ that
Bar-Natan's theory recovers Khovanov's theory by adding the
relation $t=0$,
and hence we have $$BN(F) |_{t=0} = KJ(F).$$
We remark here that Bar-Natan's theory also recovers 
Lee's theory \cite{Lee} by adding the relation $t=1$.

\section{Proof}\label{sec-pf}

For a surface-knot $F$,
taking an arbitrary point $p$ of $F$
and cutting off a small neighborhood of $p$
which is homeomorphic to the standard disk pair
$(D^4,D^2)$,
we obtain a link cobordism between an empty link and
a trivial knot.
Then we can define the following two maps 
$$\psi_F^{(
\begin{minipage}{5pt}\begin{picture}(5,5)
 \put(2.5,2.5){\circle{5}}
\end{picture}\end{minipage}
\rightarrow \emptyset)} : 
H(\bigcirc ,\mathcal{F}') \rightarrow 
H(\emptyset ,\mathcal{F}') 
\ \ \text{and} \ \ 
\psi_F^{( \emptyset \rightarrow
\begin{minipage}{5pt}\begin{picture}(5,5)
 \put(2.5,2.5){\circle{5}}
\end{picture}\end{minipage})} : 
H(\emptyset ,\mathcal{F}')
\rightarrow H(\bigcirc ,\mathcal{F}') ,$$
where $\bigcirc$ stands for a trivial knot
and the cohomology group $H(\bigcirc ,\mathcal{F}')$
of a trivial knot is $V'$, and
these two maps satisfy 
$$\psi_F = \psi_F^{(
\begin{minipage}{5pt}\begin{picture}(5,5)
 \put(2.5,2.5){\circle{5}}
\end{picture}\end{minipage}
\rightarrow \emptyset)} \circ \iota '
=\epsilon ' \circ 
\psi_F^{( \emptyset \rightarrow
\begin{minipage}{5pt}\begin{picture}(5,5)
 \put(2.5,2.5){\circle{5}}
\end{picture}\end{minipage})} .$$
For the connected sum $F_1 \# F_2$
of two surface-knots $F_1$ and $F_2$,
the map $\psi_{F_1 \# F_2}$ can be decomposed into
the composite of two maps such that
$$\psi_{F_1 \# F_2}=
\psi_{F_2}^{(
\begin{minipage}{5pt}\begin{picture}(5,5)
 \put(2.5,2.5){\circle{5}}
\end{picture}\end{minipage}
\rightarrow \emptyset)} \circ
\psi_{F_1}^{( \emptyset \rightarrow
\begin{minipage}{5pt}\begin{picture}(5,5)
 \put(2.5,2.5){\circle{5}}
\end{picture}\end{minipage})}.$$

The following two lemmas are direct consequences
of the fact that 
$$\left( m' \circ \Delta ' \right) (\p) = 2 \m  \ \,
\text{and} \ \,
\left( m' \circ \Delta ' \right) (\m) = 2t \p .$$
We note that the map $m' \circ \Delta '$ corresponds to
a link cobordism between trivial knots
induced by a trivial $T^2$-knot with two holes.

\begin{lem}\label{lem-odd}
If the surface-knot $F$ of genus $2m+1$ $(m\geq 0)$
is trivial, then we have $BN(F)=2(4t)^m$.
\end{lem}

\begin{lem}\label{lem-even}
If the surface-knot $F$ of genus $2m$ $(m\geq 0)$
is trivial, then we have 
$$\psi_F^{(
\begin{minipage}{5pt}\begin{picture}(5,5)
 \put(2.5,2.5){\circle{5}}
\end{picture}\end{minipage}
\rightarrow \emptyset)} (\m) = \pm (4t)^m.$$
\end{lem}

\begin{proof}[Proof of Theorem$~\ref{th-main}$]
Since the map $\psi_F$ induced by a surface-knot $F$
is a graded map of degree $\chi(F)$
and the degree of $t$ is $-4$,
it is easy to see the following:
\begin{itemize}
\item
If the genus of a surface-knot $F$ is $2m$ $(m \geq 0)$,
then we have $BN(F)=0$.
\item
If the genus of a surface-knot $F$ is $2m+1$ $(m \geq 0)$,
then there exists some nonnegative integer $a$
such that $BN(F)=a t^m$. 
\end{itemize}
It is sufficient to prove that the above
integer $a$ is equal to $2^{2m+1}$ for any surface-knot $F$
of genus $2m+1$.

It follows from $BN(F)= a t^m $ that
$$\psi_F^{( \emptyset \rightarrow
\begin{minipage}{5pt}\begin{picture}(5,5)
 \put(2.5,2.5){\circle{5}}
\end{picture}\end{minipage})} (1) = \pm a t^m \m .$$
Let $\Sigma_{g}$ denote a trivial surface-knot
of genus $g$.
We consider the connected sum $F \# \Sigma_{2m'}$
of $F$ and $\Sigma_{2m'}$ for a nonnegative integer $m'$.
By lemma~\ref{lem-even}, we have
$$\psi_{F \# \Sigma_{2m'}} (1) = 
\left(
\psi_{\Sigma_{2m'}}^{(
\begin{minipage}{5pt}\begin{picture}(5,5)
 \put(2.5,2.5){\circle{5}}
\end{picture}\end{minipage}
\rightarrow \emptyset)} \circ
\psi_F^{( \emptyset \rightarrow
\begin{minipage}{5pt}\begin{picture}(5,5)
 \put(2.5,2.5){\circle{5}}
\end{picture}\end{minipage})}
\right) (1)=
\pm a t^m (4t)^{m'},$$ 
and hence we have $BN(F \# \Sigma_{2m'})= a t^m (4t)^{m'}$. 

If we take the integer $m'$ such that $2m'$ is greater than 
the unknotting number \cite{HK} of $F$,
then the surface-knot $F \# \Sigma_{2m'}$ is ribbon-move 
equivalent to $\Sigma_{2(m+m')+1}$. 
When two surface-knots are related by ribbon-moves,
it is known that the induced maps on the cohomology groups 
are the same (cf. \cite{CSS,BN2}).
Hence we have  $BN(F \# \Sigma_{2m'})= 2(4t)^{m+m'}$
by Lemma~\ref{lem-odd}. This implies $a=2^{2m+1}$.
\end{proof}

\section*{Acknowledgments}

The author would like to express his sincere gratitude
to Yukio Matsumoto for encouraging him.
He would also like to thank Seiichi Kamada, Masahico Saito,
Shin Satoh for helpful comments, 
Magnus Jacobsson for telling me that Jacob Rasmussen 
has a different proof of Corollary~\ref{cor-KJ},
Isao Hasegawa and Yoshikazu Yamaguchi for stimulating discussions on 
an early version of this study.



\end{document}